\documentclass[10pt]{article}

\usepackage[utf8]{inputenc}
\usepackage{amsfonts}
\usepackage{amsmath}
\usepackage{amssymb}
\usepackage{amsthm}
\usepackage{mathrsfs}
\usepackage[american]{babel}
\usepackage[T1]{fontenc}
\usepackage{graphicx}
\usepackage{fullpage}
\usepackage{enumitem}
\usepackage{comment}

\newtheorem{theorem}{Theorem}[section]
\newtheorem{claim}[theorem]{Claim}
\newtheorem{conj}[theorem]{Conjecture}

\newtheorem{defn}[theorem]{Definition}

\def\i{{\iota}}
\def\s{{\sigma}}

\def\cA{{\cal A}}
\def\cAp{{\cA^\pr}}
\def\cB{{\cal B}}
\def\cBp{{\cB^\pr}}
\def\cIr{{\cal C}}

\def\cE{{\cal E}}
\def\cF{{\cal F}}
\def\cG{{\cal G}}
\def\cH{{\cal H}}
\def\cI{{\cal I}}
\def\cIo{{\cI^1}}
\def\cIt{{\cI^2}}
\def\cIr{{\cI^3}}
\def\cIro{{\cI^3_1}}

\def\cR{{\cal R}}
\def\cS{{\cal S}}
\def\cX{{\cal X}}
\def\cZ{{\cal Z}}

\def\oi{{\bar{i}}}
\def\oJ{{\bar{J}}}
\def\oone{{\bar{1}}}
\def\otwo{{\bar{2}}}
\def\othree{{\bar{3}}}

\def\mt{{\emptyset}}
\def\pr{{\prime}}
\def\sse{{\ \subseteq\ }}

\def\Pf{{\noindent\bf Proof.\ \ }}
\def\pf{{\noindent\it Proof.\ \ }}
\def\bpf{{\hfill $\Box$\bigskip}}
\def\dpf{{\hfill $\diamondsuit$}}

\author{
Eva Czabarka\thanks{
Department of Mathematics,
University of S. Carolina,
\texttt{czabarka@math.sc.edu}
}\\
Glenn Hurlbert\thanks{
Department of Mathematics and Applied Mathematics,
Virginia Commonwealth University
}
\thanks{
\texttt{ghurlbert@vcu.edu}.
Research partially supported by Simons Foundation Grant \#246436.
}\\
Vikram Kamat\footnotemark[2]
\ \thanks{
\texttt{vkamat@vcu.edu}
}
}

\title{Chv\'atal's conjecture for downsets of small rank}
\begin{document}
\maketitle

\begin{abstract}
A starting point in the investigation of intersecting systems of subsets of a finite set is the elementary observation that the size of a family of pairwise intersecting subsets of a finite set $[n]=\{1,\ldots,n\}$, denoted by $2^{[n]}$, is at most $2^{n-1}$, with one of the extremal structures being family comprised of all subsets of $[n]$ containing a fixed element, called as a \emph{star}.
A longstanding conjecture of Chv\'atal aims to generalize this simple observation for all \emph{downsets} of $2^{[n]}$.
In this note, we prove this conjecture for all downsets where every subset contains at most $3$ elements. 
\end{abstract}


\section{Introduction}
Let $[n]=\{1,\ldots,n\}$ and let $2^{[n]}$ (resp. $\binom{[n]}{k}$) denote the family of all subsets (resp. $r$-sized subsets) of $[n]$.
A set system containing sets of size $r$ ($r\geq 1$) is called $r$-\emph{uniform}.
Additionally, let $\binom{[n]}{\leq r}$ be the family of all subsets of size at most $r$, for any $1\leq r\leq n$.
For a family of subsets $\cF\sse 2^{[n]}$, call $\cF$ a \emph{downset} if $A\in \cF$ and $B\sse A$ implies $B\in \cF$.
Denote by $\cF^r$ those sets of $\cF$ having size $r$.
A family $\cF\sse 2^{[n]}$ is called \emph{intersecting} if $A\cap B\neq \mt$ for every $A,B\in \cF$.
For any $\cF\sse 2^{[n]}$, let $\cF_x=\{A\in \cF:x\in A\}$, called the $\cF$-star centered at $x$.
Call any $\cG\sse \cF_x$ a partial $\cF$-star centered at $x$, and call $x$ a {\it center} of such a family.
As a family may have more than one center, we call the set of all centers of $\cG$ the {\it head} of $\cG$ --- it equals the intersection of all the sets of $\cG$.
\\

A starting point in the study of intersecting set systems states that any intersecting set system on $[n]$ can contain at most $2^{n-1}$ subsets, as for any pair $(A,[n]\setminus A)$, where $A\sse [n]$, at most one can be in the intersecting family (see \cite{Ande}).
It is clear that the \emph{star} is one of the structures that attains this maximum size.
The seminal Erd\H{o}s--Ko--Rado theorem \cite{ErKoRa} proves a similar, more non-trivial result for \emph{uniform} set systems.

\begin{theorem}\label{ekr}
\cite{ErKoRa}
Let $r\leq n/2$ and let $\cF\sse \binom{[n]}{r}$ be intersecting.
Then $|\cF|\leq \binom{n-1}{r-1}$.
Furthermore, if $r<n/2$, then equality holds if and only if $\cF=\binom{n}{r}_x$, for some $x\in [n]$. 
\end{theorem}

In this note, we consider a famous longstanding conjecture of Chv\'atal (see \cite{Chva}), which deals with the ``Erd\H{o}s--Ko--Rado'' property of downsets.
Before we state the conjecture, we formulate the following definitions.
For $\cF\sse 2^{[n]}$ we set $\i(\cF)$ to be the size of the largest intersecting subfamily of $\cF$ and $\s(\cF)=\emph{max}_{x\in [n]}|\cF_x|$.
\begin{defn}[The EKR property]
A set system $\cF\sse 2^{[n]}$ is \emph{EKR} if $\i(\cF)=\s(\cF)$.
Moreover, $\cF$ is \emph{strictly} EKR if all of the largest intersecting subfamilies of $\cF$ are $\cF$-stars. 
\end{defn}

\begin{conj}\label{chvatal}
\cite{Chva}
If $\cH\sse 2^{[n]}$ is a downset, then $\cH$ is EKR. 
\end{conj}

There have been a handful of results confirming this conjecture.
For example, the trivial case $\cH=2^{[n]}$ is mentioned in \cite{Ande}, and Theorem \ref{ekr} implies the case for which $\cH=\binom{[n]}{\le k}$.
Schonheim \cite{Scho} solved the case for which the maximal elements of $\cH$ share a common element, while Chv\'atal \cite{Chva} handled the case for which the maximal sets of $\cH$ can be partitioned into two sunflowers (see definition below), each with core size 1.
In \cite{Chva} is also found the case for compressed $\cH$; Snevily \cite{Snev} strengthened this to $\cH$ being merely compressed with respect to some element (which also implies \cite{Scho}).
Miklos \cite{Mikl} (and later Wang \cite{Wang}) verified the conjecture for $\cH$ satisfying $\i(\cH)\ge |\cH|/2$, and Stein \cite{Stei} verified it for those $\cH$ having $m$ maximal sets, every $m-1$ of which form a sunflower.
Most recently, Borg \cite{Borg} solved a weighted generalization of \cite{Snev}.
\\

In this paper, we prove Conjecture \ref{chvatal} for $\cH\sse \binom{[n]}{\leq 3}$.
We also prove a slightly weaker result, one that makes an additional assumption on the size of the maximum intersecting family in $\cH$.
The advantage of this assumption is that the proof becomes significantly simpler, and the technique, which employs the famous Sunflower Lemma of Erd\H{o}s and Rado, could potentially be extended for downsets containing larger subsets. 


\subsection*{Main Results}
We verify Conjecture \ref{chvatal} for all downsets consisting of sets of size at most $3$.

\begin{theorem}\label{completechvatal}
Let $\cH\subseteq \binom{[n]}{\le 3}$ be a downset. 
Then $\cH$ is EKR. 
Moreover $\cH$ is strictly EKR, unless one of the following holds.
\begin{enumerate}\renewcommand{\labelenumi}{{\rm(\arabic{enumi})}}
\item\label{case:1} 
There is a subset $K\in\binom{[n]}{4}$ such that 
\begin{itemize}
\item
$\binom{K}{3}\subseteq \cH$, 
\item
for all $H\in\cH $, $H\subseteq K$ or $K\cap H=\emptyset$, and 
\item
the largest star in $\cH$ has size $7$.
\end{itemize}
\item\label{case:2} 
There are subsets $K\in\binom{[n]}{3}$ and (possibly empty) $M\sse [n]\setminus K$, and a subfamily $\cZ=\binom{K}{2}\cup\{Z\in\binom{K\cup M}{3}\mid |Z\cap K|=2 \}\sse\cH$ such that either
\begin{itemize}
\item 
$K\notin\cH$ and the largest star in $\cH$ has size $|\cZ|=3|M|+3$, or
\item 
$K\in\cH$ and the largest star in $\cH$ has size $|\cZ|+1=3|M|+4$.
\end{itemize}
\end{enumerate}
\end{theorem}

We also prove the following weaker result, which is significantly stronger than the result of \cite{Mikl} for subfamilies of $\binom{[n]}{\leq 3}$.

\begin{theorem}\label{bigchvatal}
Let $\cH\sse \binom{[n]}{\leq 3}$ be a downset, and let $\cI\sse \cH$ be a maximum intersecting family. If $|\cI|\geq 31$, then $\cI$ is a star.
Hence $\cH$ is EKR when $\i(\cH)\ge 31$.
\end{theorem}

Of course, some intersecting family (in particular, some star) will be so large if $|\cH|>15n$ or $|\cH^3|>10n$, for example.
\\

Our proofs use the notion of \emph{Sunflowers}, including the famous \emph{Sunflower Lemma} of Erd\H{o}s and Rado \cite{ErdRad}, as well as a variant by H{\aa}stad, et al \cite{HaJuPu}. 
We state both the Sunflower Lemma and the variant below, after the following definitions. 

\begin{defn}[Covering Set and Covering Number]
A set $S$ is a \emph{covering set} for a set system $\cF$ if $S\cap F\neq \mt$ for every $F\in \cF$.
The covering number of $\cF$, denoted by $\tau(\cF)$, is the size of the smallest covering set of $\cF$. 
\end{defn}
 
\begin{defn}[Sunflower]\label{sunflower}
A \emph{sunflower} with $k$ petals and core $C$ is a set system $\{S_1,\ldots,S_k\}$ such that for any $i\neq j$, $S_i\cap S_j=C$.
The sets $S_i\setminus C$ are the petals of the sunflower, and must be non-empty. 
If $k=1$ then we may choose $C$ to be any proper subset of $S_1$.
 \end{defn}
 
For a set system $\cF$ and set $Y$, let $\cF_Y=\{F\setminus Y:F\in \cF,Y\sse F\}$. 
\begin{defn}[$k$-Flower]\label{flower}
A $k$-flower with core $C$ is a set system $\cF$ with $\tau(\cF_C)\geq k$. 
\end{defn}
 
\begin{theorem}\label{sflemma}
\cite{ErdRad}
If a family of sets $\cF$ is $r$-uniform and $|\cF|> r!(k-1)^r$ sets, then it contains a sunflower with $k$ petals.
\end{theorem}

We will use the following variant of Theorem \ref{sflemma}.
\begin{theorem}\label{kflemma}
\cite{HaJuPu}
If $\cF$ is $r$-uniform and $|\cF|>(k-1)^r$, then $\cF$ contains a $k$-flower. 
\end{theorem}


\section{Proof of Theorem \ref{completechvatal}}

\Pf
Let $\cI$ be an intersecting subfamily of $\cH$ of maximum size. 
Our goal is to show that either $\cI$ must be a star or otherwise that $\cH$ contains a star of size equal to that of $\cI$, and to characterize the cases for which the latter happens.
\\

If $\cH$ does not contain a set of size $3$ then $\cI$ is a star unless $|\cI|=3$ and $\cI=\binom{K}{2}$ for some $K=\{x,y,z\}$. 
But then $\{\{x\},\{x,y\},\{x,z\}\}\subseteq \cH_x$, and so $|\cI|=|\cH_x|$, which is case (\ref{case:2}) of the theorem with $M=\emptyset$.
\\

Thus we may assume that $\cH$ contains a set of size $3$ and, consequently, also contains a star of size $4$. 
Therefore $|\cI|\ge 4$. 
If $\cIo\not=\emptyset$ or $|\cIt|\ge 4$ then $\cI$ is a star and we are done; so we will assume that $\cIo=\emptyset$ and $|\cIt|\le 3$ (thus $\cIr\ne\emptyset$).
Our proof splits into cases, based on $|\cIt|$.
\\

We first introduce some notation that we make use of below. 
Without loss of generality  $\bigcup_{I\in\cI^2}I=[m]$ for some $m\le 4$.
For $\emptyset\not=J\subset [m]$ define $\oJ=[m]\setminus J$, $\cA(J)=\{I\in\cIr\mid I\cap [m]=J\}$, and $C(J)=(\bigcup_{A\in\cA(J)}A)\setminus J$.
In practice, we relax the notation somewhat to write $\cA(2,3)$ instead of $\cA(\{2,3\})$, and $C(\otwo)$ instead of $C(\overline{\{2\}})$, for example.
Note that, when $m=3$, $|C(\oi)|=|\cA(\oi)|$ and $\cIr\setminus\bigcup_{i\in [3]}\cA(\oi)\subseteq\{[3]\}$.
\\

\subsection{$|\cIt|=3$}


\subsubsection{$\cIt$ is a star}

We may assume that $\cIt=\{\{1,2\},\{1,3\},\{1,4\}\}$. 
If $\cIr=\cIro$, then $\cI$ is a star.	
Otherwise, we must have $\cIr\setminus\cIro=\{\{2,3,4\}\}$. 
Therefore $(\cI\setminus\{\{2,3,4\}\})\cup\{\{1\}\}\cup\{\{1,j\}\mid j\in I\in\cIro\}$ is a star subfamily of $\cH$ that has size at least $|\cI|$ and, in fact, is larger unless $I\subseteq [4]$ for every $I\in\cIro$. 
Therefore we must have that $|\cIro|\le 3$.
\\

If $|\cIro|<3$ then, without loss of generality, $\cIro\subseteq\{\{1,2,3\},\{1,2,4\}\}$, and then $(\cI\setminus\{\{1,3\},\{1,4\}\})\cup\{\{2\},\{2,3\},\{2,4\}\}$ is a larger intersecting subfamily of $\cH$, a contradiction.
So we are left with the case in which $|\cIro|=3$ and, consequently, $|\cI|=7$ and $\cH\supseteq\binom{[4]}{\le 3}$.
\\

If there is an $H\in\cH$ such that both $H\cap [4]\not=\emptyset$ and $H\setminus [4]\not=\emptyset$ then, by taking $h\in H\cap [4]$, we have that $\binom{[4]}{\le 3}_h\cup\{H\}$ is a star in $\cH$ of size $8>7$, a contradiction.
Hence there is no such $H$, which is case (\ref{case:1}) of the theorem.
\\


\subsubsection{$\cIt$ is a triangle}

We may assume that $\cIt= \{\{1,2\},\{1,3\},\{2,3\}\}$.
\\

Relabel, if necessary, so that $0\le|C(\oone)|\le |C(\otwo)|\le |C(\othree)|$.
Then $(\cI\setminus(\cA(\oone)\cup\{\{2,3\}\}))\cup\{\{1,s\}\mid s\in C(\othree)\}\cup\{\{1\}\}$ is a star subfamily of $\cH$ of size $|\cI|+|cA(\othree)|-|\cA(\oone)|\ge|\cI|$, and so $\cH$ is EKR, and strictly so unless $|C(\oone)|=|C(\otwo)|=|C(\othree)|$, 
which we now assume.
\\

If not all the sets $C(\oi)$ are the same then, without loss of generality say $C(\oone)\not= C(\otwo)$, and so $|C(\oone)\cup C(\otwo)|>|C(\othree)|$.
Then $(\cI\setminus(\cA(\othree)\cup\{\{1,2\}\}))\cup\{\{3,s\}\mid s\in C(\oone)\cup C(\otwo)\}\cup\{\{3\}\}$ is a star subfamily of $\cH$ of size $|\cI|+|A(\oone)\cup A(\otwo)|-|\cA(\othree)|>|\cI|$, a contradiction.
\\

Finally, if $C(\oone)=C(\otwo)=C(\othree)$ then $|\cH_1|=|\cI|$, so $\cH$ is EKR, but not strictly so, giving us case (\ref{case:2}) of the theorem.
\\


\subsection{$|\cIt|=2$}

We may assume that $\cIt=\{\{1,2\},\{1,3\}\}$.
For each $I\in\cIr$ we must have $1\in I$ or $\{2,3\}\subset I$.
If $I\in\cIr\setminus(\cA(1)\cup\cA(2,3))$, then $1\in I$ and $\{2,3\}\cap I\ne\emptyset$. 
\\

If $\cA(2,3)=\emptyset$, then $\cI$ is a star, so we assume that $\cA(2,3)\ne\emptyset$. 
It must be that $\cA(1)\ne\emptyset$, since otherwise
$\cI\cup\{\{2,3\}\}$ would be a larger intersecting subfamily of $\cH$, a contradiction.
\\

Fix an $A\in\cA(1)$; then for each $k\in C(2,3)$ we must have that $k\in A$.
Thus $|C(2,3)|\le|A\setminus\{1\}|=2\le |C(1)|$. 
Hence we have that $(\cI\setminus\cA(2,3))\cup\{\{1\}\}\cup\{\{1,i\}\mid i\in C(1)\}$ is a star of size $|\cI|+|C(1)|-|\cA(2,3)|+1>|\cI|$, a contradiction.
\\ 
 

\subsection{$|\cIt|=1$}

We may assume that $\cIt=\{\{1,2\}\}$. 
Without loss of generality, both of $\cA(1), \cA(2)$ are nonempty (otherwise $\cI$ is a star and we are done). 
If, for some $i\in\{1,2\}$, we have that $|\cA(i)|\le |\cA(1,2)|$ then $(\cI\cup\{A\setminus \{i\}\mid A\in\cA(1,2)\}\cup\{\{1,2\}\setminus\{i\}\})\setminus\cA(i)$ is a star-subfamily of $\cH$ of size larger than $\cI$, which is a contradiction. 
Thus we know that $|\cA(1,2)|<\min(|\cA(1)|,|\cA(2)|)$. 
\\

If $|\cA(1)|=|\cA(2)|=1$, then $\cA(1,2)=\emptyset$, and
$(\cI\cup\{\{1,j\}:j\in C(1)\}\cup\{\{1\}\})\setminus\cA(2)$ is a star subfamily of size larger that $\cI$, a contradiction., so we may assume without loss of generality that $|\cA(1)|\ge 2$.\\

For $i\in\{1,2\}$, set $\cA^\pr(i)=\{A\setminus\{i\}\mid A\in \cA(i)\}$; then $|\cA^\pr(i)|=|\cA(i)|$.
Clearly $\cA^\pr(1)$ and $\cA^\pr(2)$ cross-intersect. 
If, for some $i\in\{1,2\}$, $\cA^\pr(i)$ is an intersecting family then $\cI\cup\cA^\pr(i)\setminus\cA(1,2)$ is an intersecting subfamily of $\cH$ that is larger that $\cI$, a contradiction, so we have that neither $\cA^\pr(1)$ nor $\cA^\pr(2)$ is intersecting. 
Since $\cA^\pr(1)$ is not intersecting and $|\cA^\pr(1)|\ge 2$, we may assume (by relabeling, if necessary) that $\{\{3,4\},\{5,6\}\}\subseteq\cA^\pr(1)$. 
Because $\cA^\pr(2)$ cross-intersects $\cA^\pr(1)$ we have $\cA^\pr(2)\subseteq\{\{3,5\},\{3,6\},\{4,5\},\{4,6\}\}$. 
In particular, $|\cA(2)|=|\cA^\pr(2)|\le 4$ and, for each $x\in\{3,4,5,6\}$, $\{1,x\}$ is a subset of some set in $\cA(1)$. 
But then $(\cI\setminus\cA(2)) \cup \{\{1,x\}\mid x\in\{3,4,5,6\}\}\cup\{\{1\}\}$ is an intersecting subfamily of $\cH$ that is larger than $\cI$, a contradiction.
\\


\subsection{$|\cIt|=0$}

Here $\cIt=\emptyset$ and $\cI$ is an intersecting family of $3$-sets such that no 2-subset of $[n]$ is contained in every element of $\cI$ (otherwise that $2$-subset could be added to $\cI$).
\\

Let $\cS$ be the largest star in $\cI$ (clearly $|\cS|\ge 2$), and let $D$ be the head of $\cS$. 
If $\cS=\cI$ then we are done, so define $\cR=\cI\setminus\cS$ and assume that $\cR\ne\emptyset$. 
In particular, for every $R\in\cR$ we must have that $R\cap D=\emptyset$; otherwise $R$ could be added to $\cS$ to create a larger star. 
If $|D|\ge 2$ then for any $d\in D$ we have that $\cI\cup\{S\setminus\{d\}\mid S\in\cS\}$ is a larger intersecting subfamily of $\cH$ than $\cI$, a contradiction. 
Therefore $|D|=1$ and, without loss of generality, $D=\{1\}$. 
\\

Let $\cF$ be the largest sunflower in $\cS$ with core $\{1\}$. 
Since any $R\in\cR$ must intersect every $F\in\cF$, we must have that  $|\cF|\le 3$. 
If $|\cF|=1$, then $\{S\setminus\{1\}\mid S\in\cS\}$ forms an intersecting family, and from the fact that $|\cS|\ge 2$ and $D=\{1\}$, we have that $\cS=\{\{1,a,b\},\{1,a,c\},\{1,b,c\}\}$ for three different numbers $a,b,c$.
Moreover, we must have $|R\cap\{a,b,c\}|\ge 2$ for every $R\in\cR$.
This means that $\cI\cup\{\{a,b\}\}$ is a larger intersecting subfamily of $\cH$ than $\cI$, a contradiction.
Therefore $2\le |\cF|\le 3$.
\\

Let $X=\left(\bigcup_{F\in\cF}F\right)$; then $|X|=2|\cF|+1$. 
Denote $X^*=X\setminus\{1\}$.
Define $Y=\left(\bigcup_{S\in\cS}S\right)\setminus X$ and set $\cS(Y)=\{S\in\cS\mid S\cap Y\ne\emptyset\}$. 
Then we must have that, for all $y\in Y$, there is an $S\in\cS(Y)$ such that $\{1,y\}\subseteq S$ and, for all $x\in X$ (including $x=1$), there is an $F\in\cF$ such that $\{1,x\}\subseteq F$. 
If $|X\cup Y|=|X|+|Y|=2|\cF|+|Y|+1>|\cR|$, then $\cS\cup\{\{1,k\}\mid k\in X\cup Y\}$ is a star subfamily of $\cH$ of size larger than $\cI$, a contradiction.
So in the rest we assume that $|\cR|\ge 2|\cF|+|Y|+1$.
\\


\subsubsection{$|\cF|=3$} 

Without loss of generality, $\cF=\{\{1,2,3\},\{1,4,5\},\{1,6,7\}\}$.
Set $\cE$ to be the family of $3$-element subsets of $X^*$ that intersect each of $\{2,3\},\{4,5\},\{6,7\}$. 
Then $|\cE|=8$ and $\cR\subseteq\cE$. 
However, if $R\in\cR$ then $X^*\setminus R\in\cE\setminus\cR$, and so $|\cR|\le 4< 7\le 2|\cF|+|Y|+1$, a contradiction. 
\\


\subsubsection{$|\cF|=2$} 

Without loss of generality, $\cF=\{\{1,2,3\},\{1,4,5\}\}$.
We have $|\cR|\ge|Y|+5$.
\\

Define $\cS^*=\cS\setminus(\cF\cup\cS(Y))$.
Clearly, $\sum_{x\in X^*}|\cS^*_x|=2|\cS^*|$, and 
$\cS^*\subseteq\{\{1,i,j\}\in\cS\mid i\in\{2,3\},j\in\{4,5\}\}$.
Denote $\cR^*=\{R\in\cR\mid R\subseteq X^*\}$. 
Since $|\cR^*|\le 4<|Y|+5$, we know that $\cR\setminus\cR^*\ne\emptyset$.
\\

For each $x\in X^*$ we set $\hat{x}$ to be the integer and $C_x$ to be the $2$-set such that $\{\{x,\hat{x}\},C_x\}=\{\{2,3\},\{4,5\}\}$ (so, in particular, $C_x=C_{\hat{x}}$).
Also define $Y_x=\{y\in Y\mid \{1,x,y\}\in\cS\}$.  
For $i\in\{2,3\}$ and $j\in\{4,5\}$ let $\cR(i,j)=\{R\in\cR\mid R\cap X^*=\{i,j\}\}$ and let $R(i,j)=\{y\mid \{i,j,y\}\in\cR(i,j)\}$. 
Note that $R(i,j)\sse Y$.
The following properties are easy to see.
\\

\begin{enumerate}
\renewcommand{\labelenumi}{\textbf{\theenumi}}
\renewcommand{\theenumi}{P.\arabic{enumi}}
\item
\label{prop:partition}
The collection $\{\cR(i,j)\mid i\in\{2,3\},j\in\{4,5\}\}$ partitions $\cR\setminus\cR^*$; in particular, at least one of these sets is nonempty. 
\item
\label{prop:sstar} 
If $\{1,\hat{i},\hat{j}\}\in\cS^*$ then $\cR(i,j)=\emptyset$. 
(Since no element of $\cR(i,j)$ intersects $\{1,\hat{i},\hat{j}\}$.)
\item
\label{prop:rrstar} 
$|\cS^*|\le 3$.
(This follows from \ref{prop:partition} and \ref{prop:sstar}.)
\item
\label{prop:qij} 
If $\min(|R(i,j)|,|R(\hat{i},\hat{j})|)\ge 1$ then $R(i,j)=R(\hat{i},\hat{j})$ with $|R(i,j)|=1$. 
Therefore if $\min(|\cR(i,j)|,$ $|\cR(\hat{i},\hat{j})|)\ge 1$ then $|\cR(i,j)|=|\cR(\hat{i},\hat{j})|)=1$.
(Since elements of $\cR(i,j)$ and $\cR(\hat{i},\hat{j}) $ can intersect in at most one element.)
\item
\label{prop:bi} 
If $X^*\setminus\{x\}\in\cR^*$ then $Y_x=\emptyset$. 
(Since $X^*\setminus\{x\}$ does not intersect sets of the form $\{1,x,y\}$ for $y\in Y$.)
\item
\label{prop:ysmall}
If $y\in Y_x$ then, for $j\in C_x$, we have $\cR(\hat{x},j)\subseteq\{\{\hat{x},j,y\}\}$.
(Since $\{1,x,y\}\in\cS$.)
\item
\label{prop:ylarge}
If $|Y_x|\ge 2$ then $\bigcup_{j\in C_x} \cR(\hat{x},j)=\emptyset$. 
(This follows from \ref{prop:ysmall}.)
\item
\label{prop:ymin} 
Since $\cR^*\ne\cR$, we have $\min(|Y_x|,|Y_{\hat{x}}|)\le 1$ for every $x\in X^*$.
(This follows from \ref{prop:partition} and \ref{prop:ylarge}.)
\item
\label{prop:last} 
If $|\cS^*|=3$ then $\cS_x^*\ne\emptyset$ for all $x\in X^*$.
\end{enumerate}

If, for some $x\in X^*$, we have $|Y_x|\ge 2$, and $|Y_{\hat{x}}|= 1$, then (from \ref{prop:partition}, \ref{prop:ysmall}, and \ref{prop:ylarge}) $|\cR\setminus\cR^*|\le 2$ and (from \ref{prop:bi}) $|\cR^*|\le 2$. 
But this means that $|\cR|\le 4< 5+|Y|$, a contradiction. 
Therefore we know that if $|Y_x|\ge 2$ then $Y_{\hat{x}}=\emptyset$.
\\

If, for some $x\in X^*$, we have $|Y_x|=|Y_{\hat{x}}|=1$ (we may assume by relabeling, if necessary, that $x=2$, so $\hat{x}=3$), then (from \ref{prop:partition} and \ref{prop:ysmall}) $|\cR\setminus\cR^*|\le 4$ and (from \ref{prop:bi}) $|\cR^*|\le 2$.
Therefore, from $Y\ne\emptyset$, we get $|\cR|\le 6\le 5+|Y|$, therefore $|\cR|=5+|Y|$ and $|Y|=1$. Without loss  of generality $Y_2=Y_3=Y=\{6\}$.
Also,  $|\cR^*|=2$, and $\cR^*=\{\{2,3,5\},\{2,3,4\}\}$ and (from \ref{prop:bi}) $Y_4=Y_5=\emptyset$; consequently $\cS(Y)=\{\{1,2,6\},\{1,3,6\}\}$.
Moreover (using \ref{prop:ysmall}), from $|\cR\setminus\cR^*|=4$ we get that, for each $i\in\{2,3\}$ and $j\in\{4,5\}$, we have $\cR(i,j)=\{\{i,j,6\}\}$.
Thus (from \ref{prop:sstar}) $\cS^*=\emptyset$. 
But this yields $|\cI_2|=5>4=|\cS|$, a contradiction.
\\

Therefore we can now assume, for all $x\in X^*$, that $\min(|Y_x|,|Y_{\hat{x}}|)=0$.
Set $L=\{x\in X^*\mid Y_x\ne\emptyset\}$. 
Then we have that $|L|\le 1$ or $L=\{i,j\}$ for some $i\in\{2,3\}$ and $j\in\{4,5\}$. 
For each $x\in X^*$ we have that
\begin{equation}
\Bigl(|\cF_x|+|Y_x|+|\cS_x^*|\Bigr) + \Bigl(\sum_{j\in C_x}|\cR(x,j)|+|\cR^*_x|\Bigr)\ \le\ |\cI_x|\ \le\ |\cS|\ ,
\label{eq:individual}\end{equation}
where we have counted the sets containing 1 before those not containing 1.
Of course, $|\cF_x|=1$ and $|S^*_x|\le 2$.
By summing over $X^*$, we obtain
$$4+\sum_{x\in X^*} |Y_x|+2|\cS^*|+2\left(\sum_{i\in\{2,3\}}\sum_{j\in\{4,5\}} |{\cal R}(i,j)|\right)+3|{\cal R}^*|\le 4|\cS|\ ,$$
which simplifies to
$$4+\sum_{x\in X^*} |Y_x|+2|\cS^*|+2|\cR|+|\cR^*|\le 4|\cS|\ .$$
In particular, 
\begin{equation}
|\cR|\le 2|\cS|-2-\frac{1}{2}\sum_{x\in X^*}|Y_x|-|\cS^*|-\frac{1}{2}|\cR^*|\ .
\label{eq:general}\end{equation}
Now we consider the following three subcases, based on the size of $L$.
\\

{\noindent\bf Case}
$L=\emptyset$
\\

Then each $Y_x=Y=\emptyset$, $\cS(Y)=\emptyset$, 
$|\cS|=|\cS^*|+2$, and $|\cR|\ge 5+|Y|=5$. 
Using \ref{prop:rrstar}, equation (\ref{eq:general}) becomes 
$$|\cR|\le 2+|\cS^*|-\frac{1}{2}|\cR^*|\le 2+3-0=5.$$
Thus $|\cR|=5$, $|\cS^*|=3$, $\cR^*=\emptyset$, $|\cS|=5$, and all inequalities hold with equality in inequality (\ref{eq:individual}). 
\\

We may assume, by relabeling, if necessary, that $\cS^*=\{\{1,2,4\},\{1,2,5\},\{1,3,4\}\}$.
Then inequality (\ref{eq:individual}) implies that, for $i\in\{2,4\}$, we have $\sum_{j\in C_i}|\cR(i,j)|=2$, and for $i\in\{3,5\}$ we have $\sum_{j\in C_i}|\cR(i,j)|=3$. 
This means that $|\cR(3,5)|-|\cR(2,4)|=1$ and, therefore, $|R(3,5)|>|R(2,4)|$.
Hence (using \ref{prop:qij}), we must have that $\cR(2,4)=\emptyset$ and $|\cR(3,5)|=1$. 
This means that $|\cR(2,5)|=|\cR(3,4)|=2$, which is a contradiction by \ref{prop:qij}. 
\\

{\noindent\bf Case}
$|L|=1$
\\

Here we may assume that $L=\{2\}$. 
Then $Y=Y_2$,  $|\cS|=2+|\cS^*|+|Y|$, $\{3,4,5\}\notin\cR^*$ (from \ref{prop:bi}, and from \ref{prop:ysmall}) $|\cR(3,j)|\le 1$ for each $j\in \{4,5\}$.
Using \ref{prop:rrstar}, equation (\ref{eq:general}) becomes 
$$|\cR|\le 2+|\cS^*|+\frac{3}{2}|Y|-\frac{1}{2}|\cR^*|\le 5+\frac{3}{2}|Y|\ .$$
\\

\begin{enumerate}
\item
$|Y|=|Y_2|=1$\\

Since $5+|Y|=6\le |\cR|\le 6+\frac{1}{2}-\frac{1}{2}|\cR^*|$, we get $|\cR|=6$, $|\cS^*|=3$, and $|\cR^*|\le 1$.
From \ref{prop:sstar} we know that there are $i\in\{2,3\}$ and $j\in\{4,5\}$ such that, if $\{\ell,m\}\ne\{i,j\}$ then $\cR(\ell,m)=\emptyset$.
Since $\cR\setminus\cR^*=\cR(i,j)$, this implies that $|\cR(i,j)|\ge 5$.
Moreover, since $|\cR(3,k)|\le 1$ for each $k\in\{4,5\}$, we have $i=2$. 
Then, using inequality (\ref{eq:individual}), we obtain $7\le 1+|Y_2|+|\cS_2^*|+\sum_{k\in C_2}|\cR(2,k)|+|\cR^*_2|\le |\cS|=6$, a contradiction.
\\

\item
$|Y|=|Y_2|\ge 2$\\

Here we have $\cR(3,4)=\cR(3,5)=\emptyset$, and so (since $\{3,4,5\}\notin\cR^*$) we get that $\sum_{j\in C_2}|\cR(2,j)|+|\cR^*_2|=|\cR|$. 
Using inequality (\ref{eq:individual}) with $x=2$ yields $1+|Y|+|\cS_2^*|+|\cR|\le 2+|\cS^*|+|Y|$.
In other words, $|\cR|\le 1+|\cS^*|-|\cS_2^*|<5$, a contradiction. 
\end{enumerate}

{\noindent\bf Case}
$|L|=2$
\\

We may assume that $L=\{2,4\}$.
Then $|\cS|=|\cS^*|+2+|Y_2|+|Y_4|$ and $Y=Y_2\cup Y_4$.
From (\ref{prop:bi}) we have $\cR^*\subseteq\{ \{2,3,4\},\{2,4,5\}\}$.
We need only consider cases for which $|\cR|\ge 5+|Y|$.
\\

\begin{enumerate}
\item
$|Y_2|,|Y_4|\ge 2$\\

From (\ref{prop:ylarge}) we know that $\cR\setminus\cR^*=\cR(2,4)$.
Using inequality (\ref{eq:individual}) with $\{i,j\}=\{2,4\}$ and the fact that $5+|Y|\le|\cR|$, we get $|Y_i|+|S_i^*|+6+|Y|\le 1+|Y_i|+|\cS_i^*|+|\cR(2,4)|+|\cR^*|\le |\cS|=|\cS^*|+2+|Y_i|+|Y_j|$, which gives, for each $j\in\{2,4\}$, that $|Y|\le |\cS^*|-4+|Y_j|< |Y_j|$, a contradiction. 

\item
$|Y_2|=1$ and $|Y_4|\ge 2$ (the case $|Y_4|=1$ and $|Y_2|\ge 2$ is handled symmetrically)\\

Without loss of generality, $Y_2=\{6\}$.
From (\ref{prop:ylarge}) we have $\cR(2,5)=\cR(3,5)=\emptyset$, and from (\ref{prop:ysmall}) we know that $\cR(3,4)\subseteq\{\{3,4,6\}\}$ and $\cR\setminus\cR^*=\cR(3,4)\cup\cR(2,4)$. 
Thus, $\cR_4=\cR$.
Also, $\cS^*\subseteq\{\{1,2,4\},\{1,2,5\},\{1,3,4\},\{1,3,5\}\}$.
Set ${\cal P}=\cS^*\setminus\cS^*_4$. 
Since $\cR(2,4)\cup\cR(3,4)=\cR\setminus\cR^*\ne\emptyset$, we get from (\ref{prop:sstar}) that $|{\cal P}|\le 1$.
Thus, $\cI_4=\cI\setminus\left(\{\{1,2,6\},\{1,2,3\}\}\cup{\cal P}\right)$, so $|\cI|\le |\cI_4|+3$. 
Therefore the family $\cI_4\cup\{\{4\},\{1,4\}\}\cup\{\{4,y\}:y\in Y_4\}\}$ is a star subfamily of $\cH$ of size $|\cI_4|+2+|Y_4|\ge |\cI_4|+4> |\cI|$, a contradiction.
\\

\item
$|Y_2|=|Y_4|=1$ (so $1\le|Y|\le 2$)\\

In this case $|\cS|=4+|\cS^*|$ and $\cR\subseteq\{\{2,3,4\},\{2,4,5\}\}$
\\

\begin{enumerate}
\item
$Y_2=Y_4=Y$\\

Here we have $|Y|=1$ so, without loss of generality, $Y=\{6\}$.
Then from \ref{prop:ysmall} we learn that $\cR(2,5)\subseteq\{\{2,5,6\}\}$, $\cR(3,5)\subseteq\{\{3,5,6\}\}$ and $\cR(3,4)\subseteq\{\{3,4,6\}\}$.
Therefore $4=6-2\le |\cR|-|\cR^*| =  |\cR\setminus\cR^*|\le 3+|\cR(2,4)|$, and so $|\cR(2,4)|\ge 1$.
By \ref{prop:sstar} we know that $\{1,3,5\}\not\in\cS^*$, and thus $\cS^*\sse \cS^*_2\cup \{\{1,3,4\}\}$.
\\

\begin{enumerate}
\item
$|\cR(2,4)|>1$\\
By \ref{prop:qij} this implies that $\cR(3,5)=\emptyset$ and, hence, for $i\in\{2,4\}$, that $\sum_{k\in C_{i}}|\cR(i,k)|\ge|\cR\setminus\cR^*|-1$.
Using equation (\ref{eq:individual}) we get that $\sum_{k\in C_2}|\cR(2,k)|+|\cR^*|+2+|\cS_2^*|\le 4+|\cS^*|$.
Since $7+|\cS_2^*|\le|\cR|+1+|\cS^*_2|\le \sum_{k\in C_2}|\cR(2,k)| +|\cR^*| +2 +|\cS_2^*|$, we obtain $3+|\cS^*_2|\le |\cS^*|\le 3$. 
This implies that $S_2^*=\emptyset$ and $|S^*|=3$, contradicting \ref{prop:last}.
\\
\medskip

\item
$|\cR(2,4)|=1$\\

Now we have (from \ref{prop:qij}) that $|\cR(3,5)|\le 1$, so we know that $|\cR(i,j)|\le 1$ for every $i\in\{2,3\},j\in\{4,5\}$.
Then $|\cR\setminus\cR^*|\ge 4$ implies that $|\cR(i,j)|=1$, $|\cR\setminus\cR^*|=4$, and $\cR^*=\{X^\pr_3,X^\pr_5\}$.
Using equation (\ref{eq:individual}) with $x\in\{2,4\}$, we get that $7+|\cS_x^*|\le 4+|\cS^*|$, and so $3\le |\cS^*|-|\cS_x^*|\le 3$, a contradiction.
\\
\end{enumerate}

\item
$Y_2\ne Y_4$\\

Here we have $|Y|=2$ and, without loss of generality, $Y_2=\{6\}$ and $Y_4=\{7\}$. 
From \ref{prop:ysmall} we get that $\cR(3,j)\subseteq\{\{3,j,6\}\}$ for each $j\in \{4,5\}$ and $\cR(i,5)\subseteq\{\{i,5,7\}\}$ for each $i\in \{2,3\}$.
This implies, in particular, that $\cR(3,5)=\emptyset$.
Thus, for $i\in\{2,4\}$, we have $\sum_{k\in C_{i}}|\cR(i,k)|\ge|\cR\setminus\cR^*|-1$. 
In particular, $$7\le |\cR| \le \sum_{k\in C_2}|\cR(2,k)| + |\cR^*| +1 \le \sum_{k\in C_2}|\cR(2,k)| + |\cR^*_2| +2\ .$$
Using inequality (\ref{eq:individual}) with $x=2$ we get that $$2 +|\cS_2^*| + \sum_{k\in C_2}|\cR(2,k)| + |\cR^*_2|\le |\cS^*|+4\ .$$
Together, these imply that $3+|\cS_2^*|\le |\cS^*|\le 3$, and so $|S^*|=3$ and $|S^*_2|=0$, which contradicts \ref{prop:last}.
\\
\end{enumerate}
\end{enumerate}

This completes the proof.
\bpf


\section{Proof of Theorem \ref{bigchvatal}} 
We now proceed to a proof of Theorem \ref{bigchvatal}.\\
\\
\Pf
Let $\cI_i=\cI\cap \binom{[n]}{i}$, for $i=1,2,3$.
We can assume $\cI_1=\mt$, since otherwise, $\cI$ is a star.
Similarly, we can assume $|\cI_2|\leq 3$.
Thus, we have $|\cI_3|\geq 28$.
Since $28=(4-1)^3+1$, we can use Theorem \ref{kflemma} to conclude that $\cI_3$ contains a $4$-flower.
Let $k\geq 4$ be maximum such that $\cS$ is a $k$-flower in $\cI_3$, and let $C$ be the core of $\cS$.
As $\cI_3$ is $3$-uniform and intersecting, every subfamily $\cG\sse \cI$ has $\tau(\cG)\leq 3$, which implies that $C\neq \mt$.
Suppose first that $C=\{a\}$, and suppose $\cI$ is not a star centered at $a$. Let $A\in \cI$ be such that $a\notin A$.
Consider the family $\cS_C$.
As $\tau(\cS_C)\geq 4$, there exists some $S_1\in \cS_C$ such that $A\cap S_1=\mt$.
Consequently, if $S^\pr=S_1\cup \{a\}$, then $S^\pr\in \cI$ and $A\cap S^\pr=\mt$, a contradiction.
As a result, we may assume that $C=\{a,b\}$.
This implies that $\cS_C$ is a family of singletons.
Consequently, $\cS$ is a sunflower with at least $4$ petals.\footnote{Note that every sunflower with $k$ petals is a $k$-flower, but the converse is not always true.}
Additionally, for every $A\in \cI_3$, $A\cap \{a,b\}\neq\mt$. 
\\

Let $\cA=\{A\in \cI_3:A\cap C=\{a\}\}$, and let $\cB=\{B\in \cI_3:B\cap C=\{b\}\}$.
We have $|\cI_3|=|S|+|\cA|+|\cB|$.
Let $\cAp=\{A-\{a\}:A\in \cA\}$, and $\cBp=\{B-\{b\}:B\in \cB\}$.
If $\cAp=\mt$ or $\cBp=\mt$, we can conclude that $\cI_3$, and hence, $\cI$ is a star (centered at either $a$ or $b$), so suppose both are non-empty.
Since $\cI$ is intersecting, $\cAp$ and $\cBp$ are cross-intersecting families, i.e. for any $A\in \cAp$ and $B\in \cBp$, $A\cap B\neq \mt.$
Let $V(\cAp)$ and $V(\cBp)$ be the vertex sets of $\cAp$ and $\cBp$ respectively, and let $n(\cX)=|V(\cX)|$ for $\cX\in \{\cAp,\cBp\}$.
We first prove the following claims.

\begin{claim}\label{clm1}
If both $\cAp$ and $\cBp$ are intersecting, or $|\cAp|\geq 2$ and $|\cBp|\geq 2$, then, 
$|\cX|\leq 2+ n(\cX)$ for each $\cX\in \{\cAp,\cBp\}$.
\end{claim}

\pf
If $\cAp$ is intersecting, it is either a triangle, or a star.
In either case, the bound follows trivially.
A similar argument works for $\cBp$, so suppose, without loss of generality that $\cAp$ has two disjoint edges, say $\{xy,x^\pr y^\pr\}$. $\cBp\sse \{xy^\pr,y^\pr y,yx^\pr,x^\pr x\}$, giving the required bound for $\cBp$.
Now, if $\cB$ has two disjoint edges, we can use a similar argument for $\cAp$, so suppose $\cBp$ is intersecting.
Without loss of generality, suppose $\cBp= \{xy^\pr,y^\pr y\}$.
Then $\cAp\sse \{xy, x^\pr y^\pr\}\cup \{A\in \binom{[n]}{2}:y^\pr\in A\}$, giving the bound $|n(\cAp)|\geq |\cAp|$.
This completes the proof of the claim.
\dpf

\begin{claim}\label{clm2}
If $\cAp$ has a pair of disjoint edges, and $|\cBp|=1$, then 
$|\cAp|\leq n(\cAp)+(|S|+1)$.
\end{claim}

\pf
Let $\{xy,x^\pr y^\pr\}$ be a pair of disjoint edges in $\cAp$, and, wlog, let $\cBp=\{xx^\pr\}$.
Let $\cAp_x=\{A\in \cAp:x\in A\}$, and let $\cAp_{x^\pr}=\{A\in \cAp:x^\pr\in A\}$.
Let $X=\{v\in [n]:v\neq x^\pr, xv\in \cA_x\}$, $X^*=\{v\in [n]:v\neq x, x^\pr v\in \cA_{x^\pr}\}$ and $R=X\cap X^*$.
Now, $|\cAp|\leq 2|R|+|X\setminus R|+|X^*\setminus R|+1$, and $n(\cAp)=2+|R|+|X\setminus R|+|X^*\setminus R|$.
So, $n(\cAp)-|\cAp|\geq -(|R|+1)$.
Since $|R|\leq |S|$ (otherwise, $R$ would be a bigger sunflower with core $\{a,x\}$ (or $\{a,x^\pr\}$), contradicting the choice of $S$), we have $n(\cAp)-|\cAp|\geq -(|S|+1)$.
\dpf
\\

In the next claim, we give lower bounds on the sizes of $\cH_a$ and $\cH_b$.

\begin{claim}
\ \\
\vspace{-0.1 in}
\begin{itemize}
\item $|\cH_a|\geq 1+ (|S|+n(\cAp)+1)+(|S|+|\cAp|).$
\item $|\cH_b|\geq 1+ (|S|+n(\cBp)+1)+(|S|+|\cBp|).$
\end{itemize}
\end{claim}

\pf
We will only give the proof for $\cH_a$, as the proof for $\cH_b$ follows identically.
We know that $|\cH_a|=\sum_{i=1}^3|\cH_a^i|$, where $\cH_a^i=\cH_a \cap \binom{[n]}{i}$ for $i\in \{1,2,3\}$.
It is trivial to note that $|\cH_a^1|=1$.
Now, consider $\cH_a^2$.
First, $\{a,b\}\in \cH_a^2$.
Also, for every $\{a,b,s\}\in S$, $\{a,s\}\in \cH_a^2$, as $\cH$ is a downset.
Similarly, for every $s\in n(\cAp)$, there exists a $t\in n(\cAp)$ such that $\{a,s,t\}\in \cI_3$, and hence, $\{a,s\}\in \cH_a^2$.
Thus, $|\cH_a^2|\geq |S|+n(\cAp)+1$.
Also, it is not hard to see that $|\cH_a^3| \geq |S|+|\cAp|$.
This completes the proof of the claim.
\dpf
\\

We will now prove that either $\cH_a$ or $\cH_b$ is bigger than $\cI$, which will complete the proof of the theorem.
It will be sufficient to prove the following claim.

\begin{claim}
$|\cH_a|+|\cH_b|> 2(|\cI_3|+3).$
\end{claim}

\pf
We will consider two cases, depending on whether or not the hypothesis of Claim \ref{clm1} is true.
Suppose the hypothesis of Claim \ref{clm1} holds, so we have $n(\cX)-|\cX|\geq -2$, for $\cX\in \{\cAp,\cBp\}$.
Thus, since $|S|>3$, we have
    \begin{eqnarray*}
    |\cH_a|+|\cH_b| & \geq & 4+4|S|+|\cAp|+|\cBp|+n(\cAp)+n(\cBp) \\
    & = & (2|S|+2|\cAp|+2|\cBp|+6)+(n(\cAp)-|\cAp|)+(n(\cBp)-|\cBp|)+2|S|-2 \\
    & \geq & 2(|\cI_3|+3)+(2|S|-6) \\
    & > & 2(|\cI_3|+3).
    \end{eqnarray*} 
    
Now, assume the hypothesis of Claim \ref{clm1} is false, so, without loss of generality, suppose $\cAp$ has a pair of disjoint edges, and $|\cBp|=1$.
Clearly, $n(\cBp)-|\cBp|=1$ and we can use Claim \ref{clm2} to conclude that $n(\cAp)-|\cAp|\geq -(|S|+1)$.
Thus, we have
  \begin{eqnarray*}
    |\cH_a|+|\cH_b| & \geq & 4+4|S|+|\cAp|+|\cBp|+n(\cAp)+n(\cBp) \\
    & \geq & (2|S|+2|\cAp|+2|\cBp|+6)-(|S|+1)+1+2|S|-2 \\
    & \geq & 2(|\cI_3|+3)+|S|-2 \\
    & > & 2(|\cI_3|+3).
    \end{eqnarray*} 
\dpf\\
\\
This proves the theorem.
\bpf
\section{Acknowledgements}
The second and third authors would like to thank Dhruv Mubayi for productive discussions on approaches to proving Theorem \ref{bigchvatal}.

\end{document}